# A Godunov type scheme for a class of LWR traffic flow models with non-local flux

Jan Friedrich*, Oliver Kolb*, Simone Göttlich*

April 27, 2018

**Abstract**

We present a Godunov type numerical scheme for a class of scalar conservation laws with non-local flux arising for example in traffic flow models. The proposed scheme delivers more accurate solutions than the widely used Lax-Friedrichs type scheme. In contrast to other approaches, we consider a non-local mean velocity instead of a mean density and provide $L^\infty$ and bounded variation estimates for the sequence of approximate solutions. Together with a discrete entropy inequality, we also show the well-posedness of the considered class of scalar conservation laws. The better accuracy of the Godunov type scheme in comparison to Lax-Friedrichs is proved by a variety of numerical examples.

**AMS subject classifications:** 35L65, 65M12, 90B20

**Keywords:** scalar conservation laws, non-local flux, Godunov scheme, traffic flow models

## 1 Introduction

Over recent years, non-local conservation laws gained growing interest for a wide field of applications such as supply chains [3], sedimentation [4], conveyor belts [12], crowd motion [8] or traffic flow [5, 10]. In the latter case, the well-known Lighthill-Whitham-Richards (LWR) model [15, 16] has been extended by considering non-local velocity terms depending on the downstream traffic so that drivers adapt their velocity to the mean traffic in front, see [5, 10].

The well-posedness of special non-local flux problems has been investigated in for example [1, 2, 7, 9]. However, only a few numerical schemes have been applied so far to solve these type of equations. The most common approach are first order Lax-Friedrichs (LxF) type schemes [2, 4, 5, 7, 10], while recently second- and higher-order schemes have been introduced [6, 11]. We remark that also these higher-order methods rely on LxF type numerical flux functions, which imply the same drawbacks known from local conservation laws. Certainly, the LxF type scheme offers a powerful tool to numerically analyze non-local flux problems but typically leads to approximate solutions with strong diffusive behavior. As we are interested in a more accurate approach, we present a Godunov type scheme for a class of scalar conservation laws with non-local flux. In addition, by deriving several properties of the scheme, we prove the well-posedness for these special non-local conservation laws, which in contrast to other models [5, 7, 10] focus on a non-local mean velocity of the downstream traffic. Furthermore, the Godunov type scheme approach allows for physically reasonable solutions meaning that a maximum principle is satisfied and negative velocities as well as negative fluxes are avoided.

*University of Mannheim, Department of Mathematics, 68131 Mannheim, Germany (janfriea@mail.uni-mannheim.de, {kolb, goettlich}@uni-mannheim.de).

This work is organized as follows: In Section 2, we present the considered class of non-local conservation laws for traffic flow. Afterwards, we introduce the Godunov type scheme and derive important properties of the scheme such as $L^\infty$ and bounded variation (BV) estimates in Section 3. Those are also used to show the well-posedness of the proposed traffic model. In Section 4, we present numerical examples, which demonstrate the better accuracy of the Godunov type scheme in comparison to the widely used LxF type scheme and also provide a comparison to the model in [7].

## 2 Modeling

We briefly present an already existing traffic flow model with non-local flux originally introduced in [5, 7, 10]. Based on the modeling ideas therein, we propose an adapted model and show its well-posedness. The key difference appears in the flux function, where instead of a mean downstream density a mean downstream velocity is considered.

### 2.1 An existing model with mean downstream density

The model considered in [7] is given by a scalar conservation law of the form

$$\partial_t \rho(t,x) + \partial_x \left(g(\rho) v(w_\eta * \rho)\right) = 0, \quad x \in \mathbb{R}, \ t > 0, \tag{1}$$

where

$$w_\eta * \rho(t,x) := \int_x^{x+\eta} \rho(t,y) w_\eta(y-x) dy, \quad \eta > 0. \tag{2}$$

For initial conditions

$$\rho(0,x) = \rho_0(x) \in \mathrm{BV}(\mathbb{R}, I), \quad I = [a,b] \subseteq \mathbb{R}^+, \tag{3}$$

the existence and uniqueness of weak entropy solutions is stated in [7, Theorem 1] if the following hypotheses are satisfied:

$$\begin{aligned}
& g \in C^1(I; \mathbb{R}^+), \\
& v \in C^2(I; \mathbb{R}^+) \qquad \text{with } v' \leq 0, \\
& w_\eta \in C^1([0,\eta]; \mathbb{R}^+) \qquad \text{with } w'_\eta \leq 0, \ \int_0^\eta w_\eta(x) dx = W_0 \ \forall \eta > 0, \lim_{\eta \to \infty} w_\eta(0) = 0.
\end{aligned} \tag{H1}$$

Note that a whole family of kernel functions $w_\eta$ is considered to also analyze the limit behaviour of the model ($\eta \to 0$ and $\eta \to \infty$).

### 2.2 Model considering a mean downstream velocity

The non-local model (1) to (3) can be applied in the context of traffic flow and chooses the velocity based on a mean downstream traffic density. In contrast to this approach, it is also reasonable to assume that drivers adapt their speed based on a mean downstream velocity, anticipating the *future* space in front of them, see Figure 1.

Therefore, we consider a slightly different model compared to (1) to (3), namely

$$\partial_t \rho(t,x) + \partial_x \left(g(\rho) \left(w_\eta * v(\rho)\right)\right) = 0, \quad x \in \mathbb{R}, \ t > 0, \tag{4}$$

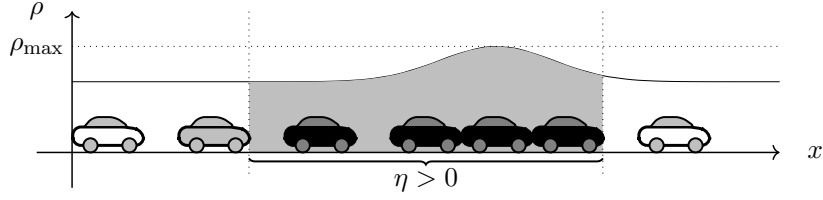

Figure 1: Illustration of a non-local traffic flow model either given by (1)-(3) or (4)-(6).

where
$$w_\eta * v(\rho)(t,x) := \int_x^{x+\eta} v(\rho(t,y))w_\eta(y-x)dy, \quad \eta > 0, \tag{5}$$

and we have given initial conditions
$$\rho(0,x) = \rho_0(x) \in \mathrm{BV}(\mathbb{R};[0,\rho_{\max}]). \tag{6}$$

For simplicity, let us also define
$$V(t,x) := w_\eta * v(\rho)(t,x). \tag{7}$$

In (4) and (5), we assume the same hypotheses as for (1) and (2) and one additional restriction:
$$\text{(H1) with } I = [0,\rho_{\max}], \quad g' \geq 0. \tag{H2}$$

As we will see in Section 3, the reformulation of the original model (1) to (3) keeps the main properties and allows for a straightforward application of a Godunov type scheme.

*Remark* 2.1. We note that in the case of a linear velocity function $v(\rho)$, e.g. $v(\rho) = 1 - \rho$, the model given by (4) and (5) coincides with (1) and (2).

The weak entropy solutions of problem (4) to (6) are intended in the following sense:

**Definition 2.2.** [13, Definition 1]
A function $\rho \in (L^1 \cap L^\infty \cap \mathrm{BV})(\mathbb{R}^+ \times \mathbb{R}; \mathbb{R})$ is a weak entropy solution if

$$\int_0^\infty \int_{-\infty}^\infty (|\rho - \kappa|\phi_t + \mathrm{sgn}(\rho - \kappa)(g(\rho) - g(\kappa))V\phi_x - \mathrm{sgn}(\rho - \kappa)g(\kappa)V_x\phi)(t,x)dxdt$$
$$+ \int_{-\infty}^\infty |\rho_0(x) - \kappa|\phi(x,0)dx \geq 0$$

for all $\phi \in C_c^1(\mathbb{R}^2;\mathbb{R}^+)$ and $\kappa \in I = [0,\rho_{\max}]$.

Our main result concerning the new model is given by the following theorem, which states the well-posedness of problem (4) to (6).

**Theorem 2.3.** *Let $\rho_0 \in \mathrm{BV}(\mathbb{R};[0,\rho_{\max}])$ and hypotheses (H2) hold. Then, the Cauchy problem*
$$\begin{cases} \partial_t \rho(t,x) + \partial_x \left(g(\rho(t,x))(v(\rho) * w_\eta)\right) = 0, & x \in \mathbb{R}, \ t > 0, \\ \rho(0,x) = \rho_0(x), & x \in \mathbb{R}, \end{cases}$$
*admits a unique weak entropy solution in the sense of Definition 2.2 and*
$$\inf_{\mathbb{R}}\{\rho_0\} \leq \rho(t,x) \leq \sup_{\mathbb{R}}\{\rho_0\} \quad \text{for a.e. } x \in \mathbb{R}, \ t > 0.$$

The proof consists of two parts: existence and uniqueness of entropy solutions. While the uniqueness proof follows from the Lipschitz continuous dependence of weak entropy solutions on the initial data, the existence proof is based on a construction of a converging sequence of approximate solutions defined by a Godunov type scheme.

## 2.3 Uniqueness of entropy solutions

One part of the proof to Theorem 2.3 is to show uniqueness of entropy solutions for the model (4) to (6). Therefore, we prove the Lipschitz continuous dependence of weak entropy solutions with respect to the initial data. Here, we follow [5, 7, 10] and use Kruzkov's doubling of variables technique [13]. Note that in the following $\|\cdot\|$ denotes $\|\cdot\|_{L^\infty}$.

**Theorem 2.4.** *Under hypotheses* (H2), *let $\rho$ and $\sigma$ be two entropy solutions of* (4) *to* (6) *with initial data $\rho_0$ and $\sigma_0$, respectively. Then, for any $T > 0$, there holds*

$$\|\rho(t,\cdot) - \sigma(t,\cdot)\|_{L^1} \leq \exp(KT)\|\rho_0 - \sigma_0\|_{L^1} \quad \forall t \in [0,T] \tag{8}$$

*with $K$ given by* (12).

*Proof.* The functions $\rho$ and $\sigma$ are weak entropy solutions of

$$\partial_t \rho(t,x) + \partial_x \left(g(\rho(t,x))V(t,x)\right) = 0, \qquad V := v(\rho) * w_\eta, \qquad \rho(0,x) = \rho_0(x),$$
$$\partial_t \sigma(t,x) + \partial_x \left(g(\sigma(t,x))U(t,x)\right) = 0, \qquad U := v(\sigma) * w_\eta, \qquad \sigma(0,x) = \sigma_0(x),$$

respectively, and $V$, $U$ are bounded measurable functions and Lipschitz continuous w.r.t. $x$ since $\rho, \sigma \in (L^1 \cap L^\infty \cap BV)(\mathbb{R}^+ \times \mathbb{R}; \mathbb{R})$. Using the classical doubling of variables technique, we get the following inequality:

$$\|\rho(t,\cdot) - \sigma(t,\cdot)\|_{L^1} \leq \|\rho_0 - \sigma_0\|_{L^1} + \|g'\| \int_0^T \int_\mathbb{R} |\rho_x(t,x)| \, |U(t,x) - V(t,x)| dx dt$$
$$+ \int_0^T \int_\mathbb{R} |g(\rho(t,x))| \, |U_x(t,x) - V_x(t,x)| dx dt, \tag{9}$$

where $\rho_x$ must be understood in the sense of measures. Applying the mean value theorem and using the properties of the kernel function, we deduce

$$|U(t,x) - V(t,x)| \leq \|v'\| \, w_\eta(0) \, \|\rho(t,\cdot) - \sigma(t,\cdot)\|_{L^1}. \tag{10}$$

Using the Leibniz integral rule and again the mean value theorem, we can also obtain for a.e. $x \in \mathbb{R}$

$$\begin{aligned}
|U_x(t,x) - V_x(t,x)| =& \Big| \int_x^{x+\eta} (v(\rho(t,y)) - v(\sigma(t,y))) w'_\eta(y-x) dy \\
&+ (v(\sigma(t,x+\eta)) - v(\rho(t,x+\eta))) w_\eta(\eta) \\
&- (v(\sigma(t,x)) - v(\rho(t,x))) w_\eta(0) \Big| \\
\leq & \|w'_\eta\| \|v'\| \|\rho(t,\cdot) - \sigma(t,\cdot)\|_{L^1} \\
&+ w_\eta(0) \|v'\| (|\rho - \sigma|(t,x+\eta) + |\rho - \sigma|(t,x)).
\end{aligned} \tag{11}$$

If we plug (10) and (11) into (9), we obtain

$$\|\rho(t,\cdot) - \sigma(t,\cdot)\|_{L^1} \leq \|\rho_0 - \sigma_0\|_{L^1} + \|v'\|\bigg(\bigg(w_\eta(0)\|g'\|\sup_{t\in[0,T]}\|\rho(t,\cdot)\|_{BV(\mathbb{R})}$$
$$+ \|w'_\eta\|\sup_{t\in[0,T]}\|g(\rho(t,\cdot))\|_{L^1}\bigg)\int_0^T\|\rho(t,\cdot) - \sigma(t,\cdot)\|_{L^1}dt$$
$$+ w_\eta(0)\sup_{t\in[0,T]}\|g(\rho(t,\cdot))\|\int_0^T\int_\mathbb{R}(|\rho - \sigma|(t,x+\eta)$$
$$+ |\rho - \sigma|(t,x))dxdt\bigg)$$
$$\leq \|\rho_0 - \sigma_0\|_{L^1} + K\int_0^T\|\rho(t,\cdot) - \sigma(t,\cdot)\|_{L^1}dt$$

with

$$K := \|v'\|\bigg(w_\eta(0)\Big(\|g'\|\sup_{t\in[0,T]}\|\rho(t,\cdot)\|_{BV(\mathbb{R})} + 2\sup_{t\in[0,T]}\|g(\rho(t,\cdot))\|\Big)$$
$$+ \|w'_\eta\|\sup_{t\in[0,T]}\|g(\rho(t,\cdot))\|_{L^1}\bigg). \quad (12)$$

By Gronwall's lemma we get (8) and for $\sigma_0 = \rho_0$ the uniqueness of entropy solutions. $\square$

The existence of weak entropy solutions is now proved in Section 3. We therefore introduce a Godunov type scheme used to construct approximate solutions.

## 3 A Godunov type scheme

The main new contribution of this work is to develop a suitable Godunov type numerical scheme for the non-local model (4) to (6). We derive $L^\infty$ and BV bounds for the approximate solutions and further provide, due to a discrete entropy inequality, all ingredients to finally prove the well-posedness result Theorem 2.3.

### 3.1 Numerical scheme

We take a space step $h$ such that for the size of the (downstream) kernel $\eta = Nh$ with $N \in \mathbb{N}$ holds. The time step is denoted by $\tau$ and we set $\lambda = \tau/h$. We define the space grid by $x_{j+\frac{1}{2}} = (j+\frac{1}{2})h$ being the cell interfaces and $x_j = jh$ the cell centers for $j \in \mathbb{Z}$, see Figure 2. The finite volume approximate solution $\bar{\rho}(t,x)$ is denoted by $\rho_j^n$ for $(t,x) \in [t^n, t^{n+1}[\times]x_{j-\frac{1}{2}}, x_{j+\frac{1}{2}}]$.

Within the proposed scheme, we intend to mimic the numerical flux function of the Godunov scheme for local conservation laws, i.e., minimizing or maximizing the flux within $[\rho_j^n, \rho_{j+1}^n]$ or $[\rho_{j+1}^n, \rho_j^n]$ depending on whether $\rho_j^n \leq \rho_{j+1}^n$ or $\rho_{j+1}^n \leq \rho_j^n$, respectively. As the convolution term $V(t,x)$ is a non-local velocity function (see (7)), there is no straightforward way of adapting this result to the non-local case. Therefore, we first examine the flux at the interface $x_{j+\frac{1}{2}}$, where the convolution term for $t \in [t^n, t^{n+1}[$ is then approximated by

$$V_{j+\frac{1}{2}}^n = \sum_{k=0}^{N-1} \gamma_k v(\rho_{j+k+1}^n) \quad (13)$$

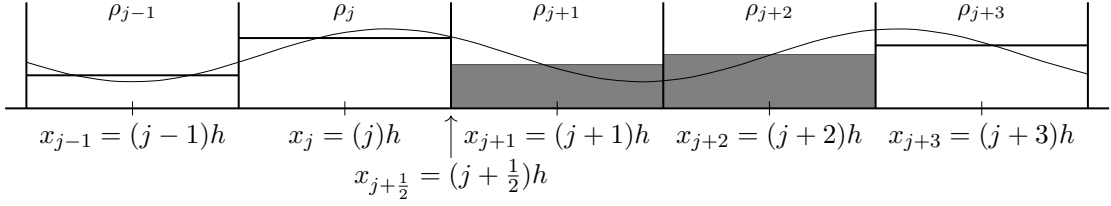

Figure 2: Space discretization and downstream kernel $\eta = Nh$ for $N = 2$ in gray.

with
$$\gamma_k = \int_{kh}^{(k+1)h} w_\eta(y) dy \qquad \forall k \in \{0, \ldots, N-1\}, \tag{14}$$
which is motivated by
$$\begin{aligned} V(t^n, x_{j+\frac{1}{2}}) &= \int_{x_{j+\frac{1}{2}}}^{x_{j+\frac{1}{2}}+\eta} w_\eta(y - x_{j+\frac{1}{2}}) v(\rho(t,y)) dy \\ &= \sum_{k=0}^{N-1} \int_{x_{j+k+\frac{1}{2}}}^{x_{j+k+\frac{3}{2}}} w_\eta(y - x_{j+\frac{1}{2}}) v(\rho(t,y)) dy \\ &\approx \sum_{k=0}^{N-1} v(\rho_{j+k+1}^n) \int_{kh}^{(k+1)h} w_\eta(y) dy \\ &= \sum_{k=0}^{N-1} \gamma_k v(\rho_{j+k+1}^n). \end{aligned} \tag{15}$$

*Remark* 3.1. For (14) we follow [4, Equation (3.2)] to satisfy
$$0 \leq V_j^n \leq v_{\max} = v(0) \quad \forall j, n$$
if $W_0 = 1$ and $\gamma_k$ is computed exactly or by an appropriate quadrature formula (such that all $\gamma_k$ are non-negative and $\sum_{k=0}^{N-1} \gamma_k = W_0 = 1$ holds).

An example for the computation of $V_{j+\frac{1}{2}}^n$ can be seen in Figure 2. The corresponding values of $\rho$, which are used in the computation of $V_{j+\frac{1}{2}}^n$ for the case $N = 2$, are gray-shaded. Based on the approximate convolution term $V_{j+\frac{1}{2}}^n$, we adapt the numerical flux function of the Godunov scheme as follows:

$$F(\rho_j^n, \ldots, \rho_{j+N}^n) = \begin{cases} \min_{\rho \in [\rho_j^n, \rho_{j+1}^n]} V_{j+\frac{1}{2}}^n g(\rho), & \text{if } \rho_j^n \leq \rho_{j+1}^n \\ \max_{\rho \in [\rho_{j+1}^n, \rho_j^n]} V_{j+\frac{1}{2}}^n g(\rho), & \text{if } \rho_j^n \geq \rho_{j+1}^n \end{cases} = V_{j+\frac{1}{2}}^n g(\rho_j^n). \tag{16}$$

Summarizing, the entire Godunov type scheme is initialized by the initial data $\rho_0$ as
$$\rho_j^0 = \frac{1}{h} \int_{x_{j-\frac{1}{2}}}^{x_{j+\frac{1}{2}}} \rho_0(x) dx \tag{17}$$

and can be computed using (13) and (16) by the finite volume scheme

$$\rho_j^{n+1} = \rho_j^n - \lambda \left( V_{j+\frac{1}{2}}^n g(\rho_j^n) - V_{j-\frac{1}{2}}^n g(\rho_{j-1}^n) \right). \tag{18}$$

*Remark* 3.2. Analogously to the scheme (17) and (18) with (13), a Godunov type scheme for the model (1) to (3) considered in [7] can be derived. For this, similar properties can be shown analogously to the following sections. One major advantage of our Godunov type scheme unlike the LxF type scheme used in [5, 7, 10] is that the numerical fluxes $F_{j+\frac{1}{2}}^n = V_{j+\frac{1}{2}}^n g(\rho_j^n)$ are always non-negative.

## 3.2 Maximum principle

The approximate solutions constructed by the Godunov type scheme (18) satisfy a strict maximum principle:

**Theorem 3.3.** *Let hypotheses* (H2) *hold. For a given initial datum $\rho_j^0$, $j \in \mathbb{Z}$ with $\rho_M = \sup_{j \in \mathbb{Z}} \rho_j^0$ and $\rho_m = \inf_{j \in \mathbb{Z}} \rho_j^0$, the approximate solutions constructed by the scheme* (18) *satisfy the bounds*

$$\rho_m \leq \rho_j^n \leq \rho_M \qquad \forall j \in \mathbb{Z},\ n \in \mathbb{N},$$

*if the following Courant-Friedrichs-Levy (CFL) condition holds:*

$$\lambda \leq \frac{1}{\gamma_0 \|v'\| \|g\| + \|v\| \|g'\|}. \tag{19}$$

*Proof.* We prove the claim per induction. For $n = 0$ the claim is obvious, so we suppose

$$\rho_m \leq \rho_j^n \leq \rho_M, \quad \forall j \in \mathbb{Z}$$

holds for a fixed $n \in \mathbb{N}$.
Before considering $\rho_j^{n+1}$ we show some general inequalities:

$$\begin{aligned} V_{j-\frac{1}{2}}^n - V_{j+\frac{1}{2}}^n &= \sum_{k=0}^{N-1} \gamma_k v(\rho_{j+k}^n) - \sum_{k=0}^{N-1} \gamma_k v(\rho_{j+1+k}^n) \\ &= \gamma_0 v(\rho_j^n) + \sum_{k=1}^{N-1} (\gamma_k - \gamma_{k-1}) v(\rho_{j+k}^n) - \gamma_{N-1} v(\rho_{j+N}) \tag{20} \\ &\leq \gamma_0 v(\rho_j^n) + \sum_{k=1}^{N-1} (\gamma_k - \gamma_{k-1}) v(\rho_M) - \gamma_{N-1} v(\rho_M) \\ &= \gamma_0 (v(\rho_j^n) - v(\rho_M)) \\ &\leq \gamma_0 \|v'\| (\rho_M - \rho_j^n), \tag{21} \end{aligned}$$

where we used the monotonicity of $w_\eta$ and $v$, and the mean value theorem.

Analogously, we obtain from (20)

$$V_{j-\frac{1}{2}}^n - V_{j+\frac{1}{2}}^n \geq \gamma_0 \|v'\| (\rho_m - \rho_j^n). \tag{22}$$

By multiplying inequality (21) by $g(\rho_M)$, subtracting $V^n_{j+\frac{1}{2}}g(\rho^n_j)$ and applying the mean value theorem, we obtain

$$V^n_{j-\frac{1}{2}}g(\rho_M) - V^n_{j+\frac{1}{2}}g(\rho^n_j) \leq \gamma_0\|v'\|\|g\|(\rho_M - \rho^n_j) + V^n_{j+\frac{1}{2}}(g(\rho_M) - g(\rho^n_j))$$
$$\leq (\gamma_0\|v'\|\|g\| + \|v\|\|g'\|)(\rho_M - \rho^n_j).$$

Therefore, under the CFL condition (19), we have

$$\rho^{n+1}_j \leq \rho^n_j + \lambda\left(V^n_{j-\frac{1}{2}}g(\rho_M) - V^n_{j+\frac{1}{2}}g(\rho^n_j)\right) \leq \rho_M.$$

Analogously, we obtain

$$V^n_{j-\frac{1}{2}}g(\rho_m) - V^n_{j+\frac{1}{2}}g(\rho^n_j) \geq (\gamma_0\|v'\|\|g\| + \|v\|\|g'\|)(\rho_m - \rho^n_j)$$

to show

$$\rho^{n+1}_j \geq \rho_m,$$

which gives us the claim. □

The maximum principle ensures that the numerical solution to (4) to (6) is bounded from above by $\rho_M = \sup_{j \in \mathbb{Z}} \rho^0_j \in [0, \rho_{\max}]$ and hence does not exceed the maximal density $\rho_{\max}$. In addition, the scheme is positivity preserving as the solution stays non-negative, since $\rho_m = \inf_{j \in \mathbb{Z}} \rho^0_j \geq 0$.

### 3.3 BV estimates

Next, we derive a BV estimate for the approximate solutions constructed by the Godunov type scheme (18). Similar to the LxF type scheme analyzed in [5, 7, 10], BV estimates cannot be derived using the standard general approaches. In particular, the Godunov type scheme also does not fit into the classical assumptions of total variation diminishing (TVD) schemes, as the total variation may slightly increase (as it is the same for the analytical solution). Nevertheless, the numerical scheme has a bounded total variation. Further, to finally also prove the existence of solutions to the model (4) to (6), we also need to provide a bound on the (discrete) variation in space and time. We begin with the BV estimate in space:

**Theorem 3.4.** *Let hypotheses* (H2) *hold,* $\rho_0 \in BV(\mathbb{R}; [0, \rho_{\max}])$ *and let* $\bar{\rho}$ *be given by* (18). *If the CFL condition* (19) *holds, then for every* $T > 0$ *the following discrete space BV estimate is satisfied:*

$$TV(\bar{\rho}(T, \cdot)) \leq \exp(C(w_\eta, v, g)T)\, TV(\rho_0)$$

*with* $C(w_\eta, v, g, \rho_{\max}) = w_\eta(0)(\|v'\|\|g\|\rho_{\max} + \|v\|\|g'\|)$.

*Proof.* Let us define

$$\Delta^n_{j+k-\frac{1}{2}} := \rho^n_{j+k} - \rho^n_{j+k-1}.$$

Then, we obtain

$$\begin{aligned}\Delta^{n+1}_{j+\frac{1}{2}} &= \Delta^n_{j+\frac{1}{2}} - \lambda\Big(V^n_{j+\frac{3}{2}}g(\rho^n_{j+1}) - 2V^n_{j+\frac{1}{2}}g(\rho^n_j) + V^n_{j-\frac{1}{2}}g(\rho^n_{j-1})\Big) \\ &= \Delta^n_{j+\frac{1}{2}} - \lambda\Big(V^n_{j+\frac{3}{2}}(g(\rho^n_{j+1}) - g(\rho^n_j)) - V^n_{j-\frac{1}{2}}(g(\rho^n_j) - g(\rho^n_{j-1})) + g(\rho^n_j)(V^n_{j+\frac{3}{2}} - 2V^n_{j+\frac{1}{2}} + V^n_{j-\frac{1}{2}})\Big) \\ &= \Delta^n_{j+\frac{1}{2}} - \lambda\Big(V^n_{j+\frac{3}{2}}g'(\xi^n_{j+\frac{1}{2}})\Delta^n_{j+\frac{1}{2}} - V^n_{j-\frac{1}{2}}g'(\xi^n_{j-\frac{1}{2}})\Delta^n_{j-\frac{1}{2}} + g(\rho^n_j)\underbrace{(V^n_{j+\frac{3}{2}} - 2V^n_{j+\frac{1}{2}} + V^n_{j-\frac{1}{2}})}_{=:(*)}\Big),\end{aligned}$$

where $\xi_{j+\frac{1}{2}}$ is between $\rho^n_j$ and $\rho^n_{j+1}$. With (20) we derive

$$\begin{aligned}(*) = &- \gamma_0 v(\rho^n_{j+1}) - \sum_{k=1}^{N-1}(\gamma_k - \gamma_{k-1})v(\rho^n_{j+1+k}) + \gamma_{N-1}v(\rho_{j+1+N}) \\ &+ \gamma_0 v(\rho^n_j) + \sum_{k=1}^{N-1}(\gamma_k - \gamma_{k-1})v(\rho^n_{j+k}) - \gamma_{N-1}v(\rho_{j+N}) \\ = &- \gamma_0 v'(\zeta_{j+\frac{1}{2}})\Delta^n_{j+\frac{1}{2}} - \sum_{k=1}^{N-1}(\gamma_k - \gamma_{k-1})v'(\zeta_{j+k+\frac{1}{2}})\Delta^n_{j+k+\frac{1}{2}} \\ &+ \gamma_{N-1}v'(\zeta_{j+N+\frac{1}{2}})\Delta^n_{j+N+\frac{1}{2}},\end{aligned}$$

where $\zeta_{j+\frac{1}{2}}$ is again between $\rho^n_j$ and $\rho^n_{j+1}$. Thus, we have

$$\begin{aligned}\Delta^{n+1}_{j+\frac{1}{2}} = &\Big(1 - \lambda\Big(V^n_{j+\frac{3}{2}}g'(\xi^n_{j+\frac{1}{2}}) - \gamma_0 v'(\zeta^n_{j+\frac{1}{2}})g(\rho^n_j)\Big)\Big)\Delta^n_{j+\frac{1}{2}} \\ &+ \lambda V^n_{j-\frac{1}{2}}g'(\xi^n_{j-\frac{1}{2}})\Delta^n_{j-\frac{1}{2}} + g(\rho^n_j)\lambda\sum_{k=1}^{N-1}(\gamma_k - \gamma_{k-1})v'(\zeta^n_{j+k+\frac{1}{2}})\Delta^n_{j+k+\frac{1}{2}} \\ &- g(\rho^n_j)\lambda v'(\zeta^n_{j+\frac{1}{2}})\gamma_{N-1}\Delta^n_{j+N+\frac{1}{2}}.\end{aligned}$$

Due to the CFL condition (19) and hypotheses (H2), all terms before the differences are positive and we get

$$\begin{aligned}\sum_j |\Delta^{n+1}_{j+\frac{1}{2}}| \leq &\sum_j \Big(1 - \lambda\Big(V^n_{j+\frac{3}{2}}g'(\xi^n_{j+\frac{1}{2}}) - \gamma_0 v'(\zeta^n_{j+\frac{1}{2}})g(\rho^n_j)\Big)\Big)|\Delta^n_{j+\frac{1}{2}}| \\ &+ \lambda\sum_j V^n_{j-\frac{1}{2}}g'(\xi^n_{j-\frac{1}{2}})|\Delta^n_{j-\frac{1}{2}}| \\ &+ \lambda\sum_j g(\rho^n_j)\sum_{k=1}^{N-1}(\gamma_k - \gamma_{k-1})v'(\zeta^n_{j+k+\frac{1}{2}})|\Delta^n_{j+k+\frac{1}{2}}| \\ &- \lambda\sum_j g(\rho^n_j)v'(\zeta^n_{j+\frac{1}{2}})\gamma_{N-1}|\Delta^n_{j+N+\frac{1}{2}}|.\end{aligned}$$

Rearranging the indices we obtain

$$\sum_j |\Delta_{j+\frac{1}{2}}^{n+1}| \leq \sum_j \bigg(1 - \lambda(V_{j+\frac{3}{2}}^n - V_{j+\frac{1}{2}}^n)g'(\xi_{j+\frac{1}{2}}^n)$$
$$- v'(\zeta_{j+\frac{1}{2}}^n)\lambda\bigg(-\gamma_0 g(\rho_j^n) + \sum_{k=1}^{N-1}(\gamma_{k-1} - \gamma_k)g(\rho_{j-k}^n) + \gamma_{N-1}g(\rho_{j-N}^n)\bigg)\bigg)|\Delta_{j+\frac{1}{2}}^n|$$
$$\leq \sum_j \bigg(1 - \lambda(V_{j+\frac{3}{2}}^n - V_{j+\frac{1}{2}}^n)g'(\xi_{j+\frac{1}{2}}^n) + \gamma_0\|v'\|\|g\|\bigg)|\Delta_{j+\frac{1}{2}}^n|.$$

Using inequality (21), for which

$$(21) \leq \gamma_0 \|v\| \rho_{\max}$$

holds, and with $\gamma_0 \leq h w_\eta(0)$, we obtain

$$\sum_j |\Delta_{j+\frac{1}{2}}^{n+1}| \leq \big(1 + \lambda\gamma_0(\|v\|\|g'\|\rho_{\max} + \|v'\|\|g\|)\big)\sum_j |\Delta_{j+\frac{1}{2}}^n|$$
$$\leq \big(1 + \tau w_\eta(0)(\|v\|\|g'\|\rho_{\max} + \|v'\|\|g\|)\big)\sum_j |\Delta_{j+\frac{1}{2}}^n|.$$

Therefore, we recover the following estimate for the total variation

$$TV(\bar\rho(T, \cdot)) \leq \big(1 + \tau w_\eta(0)(\|v\|\|g'\|\rho_{\max} + \|v'\|\|g\|)\big)^{T/\tau} TV(\bar\rho(0, \cdot))$$
$$\leq \exp\big(w_\eta(0)(\|v\|\|g'\|\rho_{\max} + \|v'\|\|g\|)T\big) TV(\rho_0). \qquad (23)$$

$\square$

We are now able to also provide an estimate for the discrete total variation in space and time:

**Theorem 3.5.** *Let hypotheses* (H2) *hold,* $\rho_0 \in BV(\mathbb{R}; [0, \rho_{\max}])$ *and let* $\bar\rho$ *be given by* (18). *If the CFL condition* (19) *holds, then for every* $T > 0$ *the following discrete space and time total variation estimate is satisfied:*

$$TV(\bar\rho; \mathbb{R} \times [0, T]) \leq T \exp\big(C(w_\eta, v, g, \rho_{\max})T\big)\big(1 + W_0\|v'\|\|g\| + \|v\|\|g'\|\big) TV(\rho_0)$$

*with* $C(w_\eta, v, g, \rho_{\max}) = w_\eta(0)(\|v'\|\|g\|\rho_{\max} + \|v\|\|g'\|)$.

*Proof.* We fix $T \in \mathbb{R}^+$. If $T \leq \tau$, then $TV(\bar\rho; \mathbb{R} \times [0, T]) \leq T \cdot TV(\rho_0)$. For $T > \tau$ let $M \in \mathbb{N} \setminus \{0\}$ such that $M\tau < T \leq (M+1)\tau$. Then

$$TV(\bar\rho; \mathbb{R} \times [0, T]) = \underbrace{\sum_{n=0}^{M-1}\sum_j \tau|\rho_{j+1}^n - \rho_j^n| + (T - M\tau)\sum_j |\rho_{j+1}^M - \rho_j^M|}_{\leq T\exp(C(w_\eta, v, g, \rho_{\max})T)TV(\rho_0)} + \sum_{n=0}^{M-1}\sum_j h|\rho_j^{n+1} - \rho_j^n|.$$

If we consider the scheme (18), we obtain

$$\begin{aligned}
\rho_j^{n+1} - \rho_j^n &= \lambda\Big(V_{j-\frac{1}{2}}^n g(\rho_{j-1}^n) - V_{j+\frac{1}{2}}^n g(\rho_j^n)\Big) \\
&= \lambda\Big((V_{j-\frac{1}{2}}^n - V_{j+\frac{1}{2}}^n)g(\rho_{j-1}^n) - V_{j+\frac{1}{2}}^n(g(\rho_j^n) - g(\rho_{j-1}^n))\Big) \\
&= \lambda\Big(-g(\rho_{j-1}^n)\sum_{k=0}^{N-1}\gamma_k v'(\zeta_{j+k+\frac{1}{2}}^n)(\rho_{j+k+1}^n - \rho_{j+k}^n) - V_{j+\frac{1}{2}}^n g'(\xi_{j+\frac{1}{2}}^n)(\rho_j^n - \rho_{j-1}^n)\Big).
\end{aligned}$$

Taking absolute values yields

$$|\rho_j^{n+1} - \rho_j^n| \leq \lambda\Big(\|v'\|\|g\|\sum_{k=0}^{N-1}\gamma_k|\rho_{j+k+1}^n - \rho_{j+k}^n| + \|v\|\|g'\|\,|\rho_j^n - \rho_{j-1}^n|\Big).$$

Summing over $j$ and rearranging the indices gives us

$$\sum_j h|\rho_j^{n+1} - \rho_j^n| \leq \tau \sum_j |\rho_j^n - \rho_{j-1}^n|\big(\|v'\|\|g\|W_0 + \|v\|\|g'\|\big)$$

so that we have

$$\sum_{n=0}^{M-1}\sum_j h|\rho_j^{n+1} - \rho_j^n| \leq T\exp\big(C(w_\eta, v, g, \rho_{\max})T\big)\big(\|v'\|\|g\|W_0 + \|v\|\|g'\|\big)TV(\rho_0).$$

Therefore, we recover

$$TV(\bar\rho; \mathbb{R}\times[0,T]) \leq T\exp\big(C(w_\eta, v, g, \rho_{\max})T\big)\big(1 + W_0\|v'\|\|g\| + \|v\|\|g'\|\big)TV(\rho_0)$$

as desired. $\square$

### 3.4 Discrete entropy inequality

As another desirable property and final ingredient regarding the proof of Theorem 2.3, we next show that the approximate solutions obtained by the Godunov type scheme (18) fulfill a discrete entropy inequality. Therefore, we follow [2, 5, 7, 10] and define

$$G_{j+\frac{1}{2}}(u) := V_{j+\frac{1}{2}}^n g(u), \qquad F_{j+\frac{1}{2}}^\kappa(u) := G_{j+\frac{1}{2}}(u \wedge \kappa) - G_{j+\frac{1}{2}}(u \vee \kappa)$$

with $a \wedge b = \max(a,b)$ and $a \vee b = \min(a,b)$.

**Theorem 3.6.** *Let $\rho_j^n$, $j \in \mathbb{Z}$, $n \in \mathbb{N}$ be given by (18), and let the CFL condition (19) and hypotheses (H2) hold. Then we have*

$$|\rho_j^{n+1} - \kappa| - |\rho_j^n - \kappa| + \lambda(F_{j+\frac{1}{2}}^\kappa(\rho_j^n) - F_{j-\frac{1}{2}}^\kappa(\rho_{j-1}^n)) \qquad (24)$$
$$+ \lambda\operatorname{sgn}(\rho_j^{n+1} - \kappa)g(\kappa)(V_{j+\frac{1}{2}}^n - V_{j-\frac{1}{2}}^n) \leq 0$$

*for all $j \in \mathbb{Z}$, $n \in \mathbb{N}$ and $\kappa \in I = [0, \rho_{\max}]$.*

*Proof.* The proof closely follows [2, 5, 7]. We set

$$\tilde{H}_j(u,w) = w - \lambda(G_{j+\frac{1}{2}}(w) - G_{j-\frac{1}{2}}(u))$$
$$= w - \lambda(V^n_{j+\frac{1}{2}} g(w) - V^n_{j-\frac{1}{2}} g(u)),$$

which is a monotone non-decreasing function with respect to each variable under the CFL condition (19) since we have

$$\frac{\partial \tilde{H}_j}{\partial w} = 1 - \lambda V^n_{j+\frac{1}{2}} g'(w) \geq 0, \quad \frac{\partial \tilde{H}_j}{\partial u} = \lambda V^n_{j-\frac{1}{2}} g'(u) \geq 0.$$

Moreover, we have the identity

$$\tilde{H}_j(\rho^n_{j-1} \wedge \kappa, \rho^n_j \wedge \kappa) - \tilde{H}_j(\rho^n_{j-1} \vee \kappa, \rho^n_j \vee \kappa)$$
$$= |\rho^n_j - \kappa| - \lambda(F^\kappa_{j+\frac{1}{2}}(\rho^n_j) - F^\kappa_{j-\frac{1}{2}}(\rho^n_{j-1})).$$

Then, by monotonicity, the definition of the scheme (18) and by using (for the last inequality) the non-negativity of $(a,b) \mapsto (\text{sgn}(a+b) - \text{sgn}(a))(a+b)$, we get (24):

$$\tilde{H}_j(\rho^n_{j-1} \wedge \kappa, \rho^n_j \wedge \kappa) - \tilde{H}_j(\rho^n_{j-1} \vee \kappa, \rho^n_j \vee \kappa)$$
$$\geq \tilde{H}_j(\rho^n_{j-1}, \rho^n_j) \wedge \tilde{H}_j(\kappa, \kappa) - \tilde{H}_j(\rho^n_{j-1}, \rho^n_j) \vee \tilde{H}_j(\kappa, \kappa)$$
$$= |\tilde{H}_j(\rho^n_{j-1}, \rho^n_j) - \tilde{H}_j(\kappa, \kappa)|$$
$$= \text{sgn}(\tilde{H}_j(\rho^n_{j-1}, \rho^n_j) - \tilde{H}_j(\kappa, \kappa)) \cdot (\tilde{H}_j(\rho^n_{j-1}, \rho^n_j) - \tilde{H}_j(\kappa, \kappa))$$
$$= \text{sgn}(\tilde{H}_j(\rho^n_{j-1}, \rho^n_j) - \kappa + \lambda g(\kappa)(V^n_{j+\frac{1}{2}} - V^n_{j-\frac{1}{2}})) \cdot (\tilde{H}_j(\rho^n_{j-1}, \rho^n_j) - \kappa + \lambda g(\kappa)(V^n_{j+\frac{1}{2}} - V^n_{j-\frac{1}{2}}))$$

$$\geq \text{sgn}(\tilde{H}_j(\rho^n_{j-1}, \rho^n_j) - \kappa) \cdot (\tilde{H}_j(\rho^n_{j-1}, \rho^n_j) - \kappa + \lambda g(\kappa)(V^n_{j+\frac{1}{2}} - V^n_{j-\frac{1}{2}}))$$
$$= |\tilde{H}_j(\rho^n_{j-1}, \rho^n_j) - \kappa| + \lambda \text{sgn}(\tilde{H}_j(\rho^n_{j-1}, \rho^n_j) - \kappa) g(\kappa)(V^n_{j+\frac{1}{2}} - V^n_{j-\frac{1}{2}})$$
$$= |\rho^{n+1}_j - \kappa| + \lambda \text{sgn}(\rho^{n+1}_j - \kappa) g(\kappa)(V^n_{j+\frac{1}{2}} - V^n_{j-\frac{1}{2}}).$$

□

## 3.5 Proof of Theorem 2.3

Since we have already shown uniqueness of weak entropy solutions to the model (4) to (6), it remains to finalize the existence proof. Similar to [5, Section 4] and [7, Proof of Theorem 1], the convergence of the approximate solutions constructed by the Godunov type scheme (18) towards the unique weak entropy solution can be proven by applying Helly's theorem. The latter can be applied due to Theorems 3.3 and 3.5 and states that there exists a sub-sequence of the constructed $\bar{\rho}$ that converges to some $\rho \in (L^1 \cap L^\infty \cap BV)(\mathbb{R}^+ \times \mathbb{R}; I)$. Following a Lax-Wendroff type argument, see [14, Theorem 12.1], one can show that the limit function $\rho$ is a weak entropy solution of (4) to (6) in the sense of Definition 2.2. Together with the uniqueness result in Theorem 2.4, this concludes the proof of Theorem 2.3.

# 4 Numerical examples

In this section, we present some numerical examples demonstrating the advantages of the Godunov type scheme in comparison to the widely used LxF type scheme. The latter will be briefly introduced in the following section. In addition, we also comment on the differences between the model considered in this work (4) to (6) and the earlier one (1) to (3).

## 4.1 A LxF type scheme

The common scheme used so far for the problem (1) to (3) is a LxF type scheme, where the downstream velocity of the convolution term is computed by

$$V_j^n = v\left(h \sum_{k=0}^{N-1} w_\eta^k \rho_{j+k}^n\right) \tag{25}$$

with $w_\eta^k = w_\eta(kh)$. This discretization of $w_\eta$ implies (for $W_0 = 1$)

$$h \sum_{k=0}^{N-1} w_\eta^k \leq 1 + w_\eta(0)h.$$

Thus, the approximation (25) of the convolution term slightly overestimates the traffic density and therewith underestimates the velocity. Further, unphysical densities beyond $\rho_{\max}$ or negative velocities are possible. This can be both avoided by discretizing the kernel function as proposed in (14). In the following subsections we will use this discretization to avoid that the accuracy studies are biased by different quadrature rule errors.

The numerical flux function of the LxF scheme is given by

$$F_{j+\frac{1}{2}}^n := \frac{V_j^n g(\rho_j^n) + V_{j+1}^n g(\rho_{j+1}^n)}{2} + \frac{\alpha}{2}(\rho_j^n - \rho_{j+1}^n)$$

with $\alpha \geq 0$ being the viscosity coefficient. This leads to the scheme

$$\rho_j^{n+1} = \rho_j^n + \frac{\lambda\alpha}{2}(\rho_{j-1}^n - 2\rho_j^n + \rho_{j+1}^n) + \frac{\lambda}{2}(V_{j-1}^n g(\rho_{j-1}^n) - V_{j+1}^n g(\rho_{j+1}^n)). \tag{26}$$

For the corresponding CFL condition and restrictions on $\alpha$, we refer to [7, Proposition 2].

*Remark* 4.1. Note that the LxF type scheme can be adapted to the model (4) and vice versa the Godunov type scheme to model (1), where in comparison the LxF type scheme adds more diffusion to the numerical solution (see also Section 4.2).

## 4.2 Accuracy of the Godunov type scheme

In the following sections, we consider the solution of model (4) to (6) and model (1) to (3) in the case of traffic flow. So we have $g(\rho) = \rho$, a velocity function $v(\rho)$ specified below, and a road of length $L = 1$. We are interested in the solution at a certain final time $T$ for the initial conditions

$$\rho_0(x) = \begin{cases} 1, & \text{if } \frac{1}{3} \leq x \leq \frac{2}{3}, \\ \frac{1}{3}, & \text{else.} \end{cases} \tag{27}$$

For simplicity, we use periodic boundary conditions in all our examples. To compute the $L^1$ error of an approximate solution $\rho_h$ with step size $h$ compared to a reference solution $\rho_{\tilde{h}}$ with step size $\tilde{h}$ at time $T$, we apply

$$L^1 \text{ error} := h \sum_j |\rho_h(T, x_j) - \rho_{\tilde{h}}(T, x_j)|.$$

Now, in order to compare the accuracy of the Godunov scheme and the LxF scheme, we first use the linear velocity function

$$v(\rho) = 1 - \rho,$$

as for this choice both models coincide (see Remark 2.1) and the different discretization schemes both have been well analyzed. We consider the final time $T = 0.1$ and the kernel function $w_\eta(x) = 3(\eta^2 - x^2)/(2\eta^3)$ with $\eta = 0.1$. We compare the $L^1$ distances of the LxF and the Godunov type scheme to a reference solution computed with the LxF type scheme and $\tilde{h} = 0.02 \cdot 2^{-9}$. The spatial step size is given by $h = 0.02 \cdot 2^{-n}$ with $n \in \{0, \ldots, 6\}$. The time step parameter $\tau$ is given by the minimum of the CFL condition (19) and the CFL condition of the LxF type scheme.

The results for this first test case are given in Table 1. Obviously, the $L^1$ errors of the Godunov type scheme are significantly smaller than the ones of the LxF scheme.

Table 1: $L^1$ errors for $v(\rho) = 1 - \rho$ at $T = 0.1$.

| $n$ | Godunov  | LxF      |
|-----|----------|----------|
| 0   | 9.38e-03 | 1.99e-02 |
| 1   | 6.97e-03 | 1.30e-02 |
| 2   | 4.29e-03 | 9.31e-03 |
| 3   | 3.00e-03 | 6.41e-03 |
| 4   | 1.96e-03 | 4.27e-03 |
| 5   | 1.33e-03 | 2.71e-03 |
| 6   | 9.05e-04 | 1.64e-03 |

The better accuracy of the Godunov type scheme can also be seen directly in Figure 3. We notice that in the presence of the two discontinuities, the Godunov type scheme in particular shows a better resolution of the solution structure than the LxF scheme at the left-hand side, while the resolution for the rarefaction wave close to the jump on the right-hand side is quite similar for the short time period $T = 0.1$.

If we consider a longer time period, i.e. $T = 1$, the good performance of the Godunov type scheme can be observed in Figure 4. Here, the reference solution is computed with the LxF type scheme with a spatial step size of $\tilde{h} = 0.02 \cdot 2^{-7}$.

In addition, we consider the accuracy for a non-linear velocity function,

$$v(\rho) = 1 - \rho^5.$$

We choose the constant kernel $w_\eta(x) = \frac{1}{\eta}$ with $\eta = 0.1$ and consider again the solution of the initial conditions (27) but at time $T = 0.05$ on a road of length $L = 1$. The spatial step size is given as above by $h = 0.02 \cdot 2^{-n}$ with $n \in \{0, \ldots, 6\}$. The reference solution is computed by the LxF scheme adapted to the model (4) and with $\tilde{h} = 0.02 \cdot 2^{-9}$.

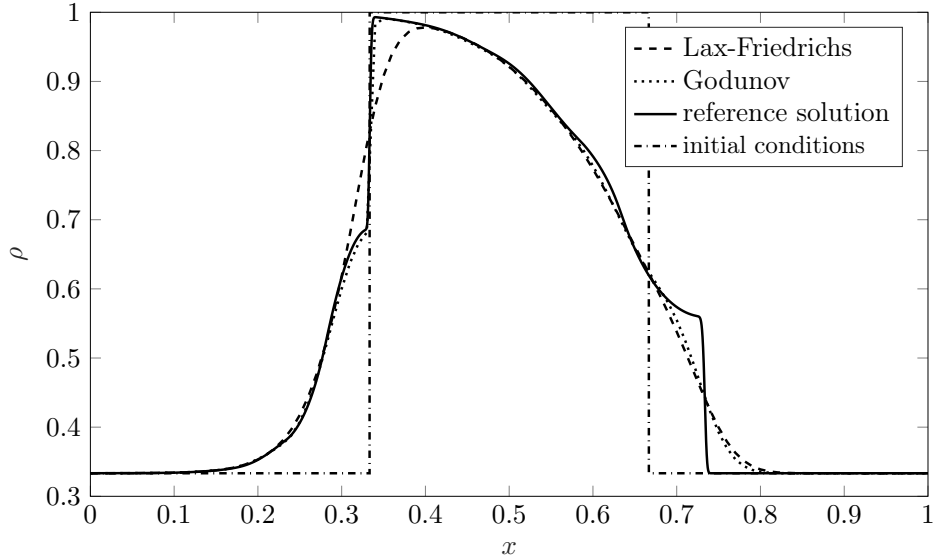

Figure 3: Comparison of the Godunov and LxF scheme for $v(\rho) = 1 - \rho$, $h = 0.01$ at $T = 0.1$.

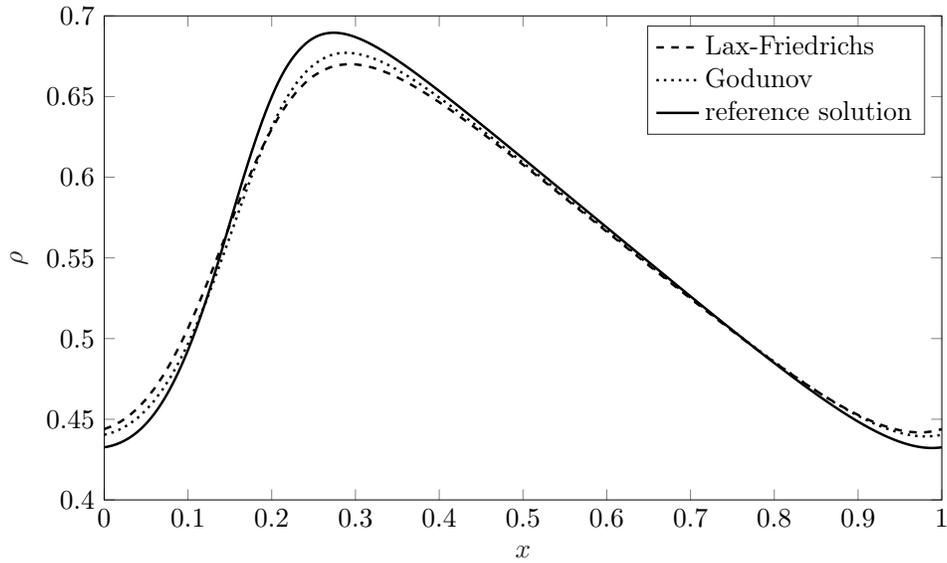

Figure 4: Comparison of the Godunov and LxF scheme for $v(\rho) = 1 - \rho$, $h = 0.01$ at $T = 1$.

The results for the non-linear test case can be seen in Table 2. Similar to the linear test case, the $L^1$ errors of the Godunov type scheme are significantly smaller than the ones of the LxF scheme. In addition, the better resolution of the Godunov type scheme can be seen in Figure 5. Again, similar to the linear case, the Godunov type scheme in particular shows a better resolution of the solution structure than the LxF scheme at the left-hand side, while the resolution for the rarefaction wave close to the jump on the right-hand side is quite similar.

Table 2: $L^1$ errors for $v(\rho) = 1 - \rho^5$ at $T = 0.05$.

| $n$ | Godunov | LxF |
|---|---|---|
| 0 | 1.77e-02 | 3.13e-02 |
| 1 | 1.24e-02 | 2.20e-02 |
| 2 | 8.49e-03 | 1.41e-02 |
| 3 | 5.18e-03 | 8.67e-03 |
| 4 | 3.29e-03 | 5.45e-03 |
| 5 | 2.02e-03 | 3.47e-03 |
| 6 | 1.21e-03 | 2.06e-03 |

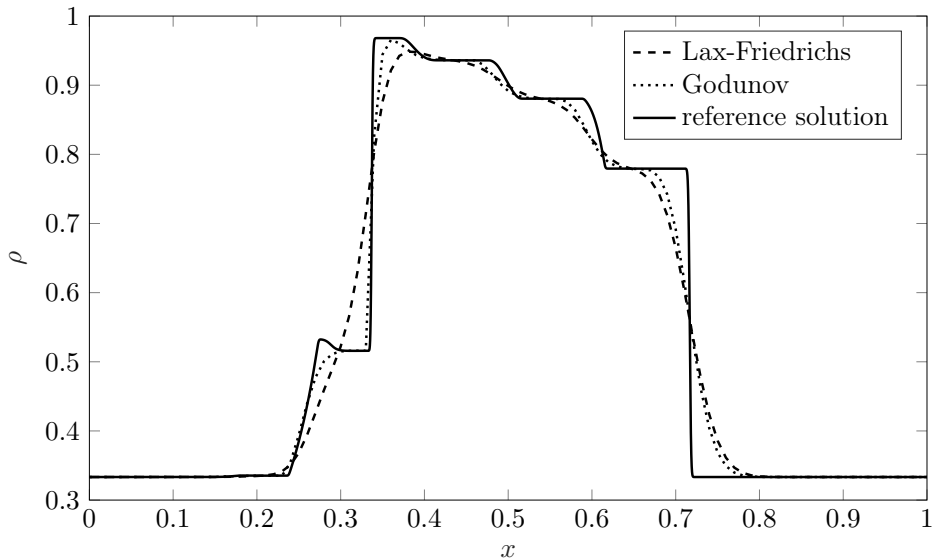

Figure 5: Comparison of the Godunov and LxF scheme for $v(\rho) = 1-\rho^5$, $h = 0.01$ at $T = 0.05$.

### 4.3 Comparison of the models

Next, we aim to discuss the differences between the models (4) to (6) and (1) to (3). To see the different dynamics within the two models, we have to choose a non-linear velocity function and we choose the same non-linear velocity $v(\rho) = 1 - \rho^5$ as before with the same parameter $\eta = 0.1$ for the constant kernel function and final time $T = 0.05$.

For a fair comparison of the evolution of densities, we apply the Godunov type scheme to both models. Figure 6 shows the results for a spatial step size $h = 0.01$ and a time step size $\tau$ determined by the CFL condition. It can be observed that the approximate solution of the earlier model (1) to (3) contains rather large oscillations while resolving the traffic jam. In contrast to that, the model (4) to (6) resolves the traffic jam in a more monotone way.

### 4.4 Limit $\eta \to 0$

Finally, we take a look at the limit case $\eta \to 0$ and investigate numerically whether the approximate solutions constructed by the proposed Godunov type scheme converge towards the solution of the (local) LWR traffic model. Note that this property is far from obvious since the constants in Theorems 2.4 and 3.5 blow up (see also [7]). For a numerical investigation, we consider the same (non-linear) scenario as above with a fixed space step size $h = 0.5 \cdot 10^{-4}$ and vary $\eta \in \{10^{-1}, 10^{-2}, 10^{-3}, 10^{-4}\}$. As final time we take $T = 0.05$.

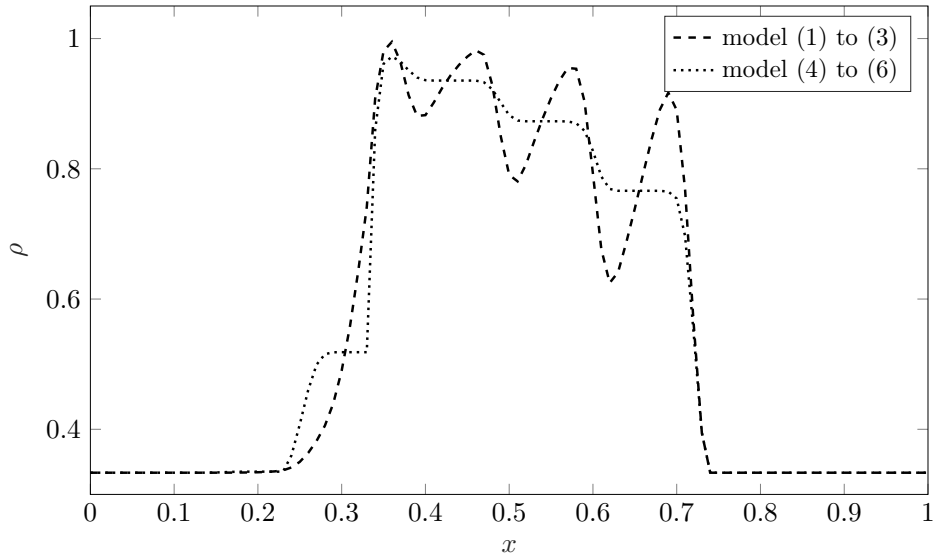

Figure 6: Approximate solutions at $T = 0.05$ for the two models with non-linear velocity function $v(\rho) = 1 - \rho^5$.

To evaluate the convergence, we compute the $L^1$ distances between the approximate solutions obtained for the proposed Godunov type scheme applied to (4) to (6) and the result of a classical Godunov scheme for the corresponding local LWR problem. Obviously, the corresponding $L^1$ distances shown in Table 3 demonstrate the convergence towards the solution of the local traffic model. The results are further illustrated in Figure 7.

Table 3: $L^1$ distances between the approximate solutions to the local LWR model and the non-local model for different $\eta$ at $T = 0.05$.

| $\eta$ | $10^{-1}$ | $10^{-2}$ | $10^{-3}$ | $10^{-4}$ |
|---|---|---|---|---|
| $L^1$ distance | 4.46e-02 | 6.85e-03 | 9.90e-04 | 1.60e-04 |

## 5 Conclusion

In this work, we have presented a Godunov type scheme for a class of non-local conservation laws. For this novel scheme we provide $L^\infty$ and BV bounds as well as a discrete entropy inequality. Based on these results, we also proved the well-posedness, i.e., existence and uniqueness of weak entropy solutions. The proposed Godunov type scheme can be adapted to other classes of non-local conservation laws and is very promising as it adds less numerical diffusion to the solution compared to the LxF type scheme. Certainly, this advantage can also be exploited by using the underlying numerical flux function within higher-order methods.

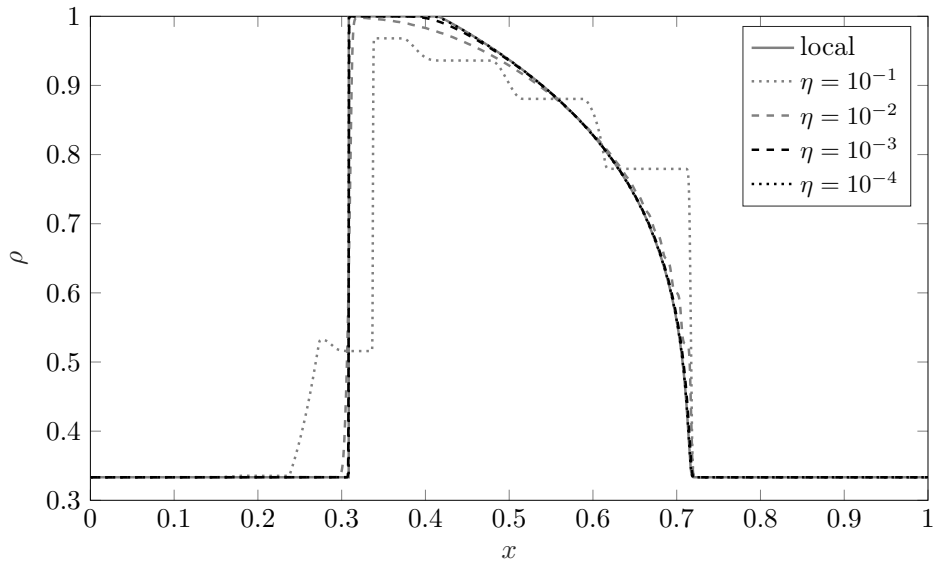

Figure 7: Approximate solutions to the LWR and non-local model (4) to (6) for different $\eta$ at $T = 0.05$.

In future work we aim at constructing several higher order methods based on the presented Godunov type scheme. In addition, the considered model with mean downstream velocity may be advantageous in the context of networks. Here we aim at investigating appropriate coupling conditions and suitable discretization schemes.